\documentclass[12pt]{article}

\usepackage{amssymb,amsfonts }
\usepackage{ytableau}
\usepackage{amssymb}
\usepackage[margin=1in]{geometry}
\usepackage{mathrsfs}
\usepackage{stmaryrd}
\usepackage{amsfonts,amsmath,amssymb,amscd}
\usepackage{shadow}
\usepackage{graphicx}
\usepackage{color}
\usepackage{longtable}
\usepackage{pstricks,multido}
\usepackage{hyperref}
\usepackage{qtree}
\usepackage{tikz}\usetikzlibrary{matrix,arrows}
\usepackage[utf8]{inputenc}
\usepackage[T1]{fontenc}
\allowdisplaybreaks

\newtheorem{thm}{Theorem}[section]
\newtheorem{lem}[thm]{Lemma}

\newtheorem{cor}[thm]{Corollary}
\newtheorem{remark}{Remark}[section]
\newtheorem{example}{Example}[section]

\newtheorem{conj}{Conjecture}[section]

\def\qed{\hfill \rule{4pt}{7pt}}

\def\pf{\noindent {\it{Proof.} \hskip 2pt}}

\numberwithin{equation}{section}

\newcommand{\PATTERN}{
    \draw[step=1, xshift=14pt, yshift=14pt, \cfill, line cap=round] (0,0) grid (3,3);
    \draw[step=1, xshift=14pt, yshift=14pt, thick] (0,1) -- (3,1);
    \draw[step=1, xshift=14pt, yshift=14pt, thick] (1,0) -- (1,3);
    \foreach \x/\y in {1/1,2/3,3/2} \node[disc, fill=black] at (\x,\y) {};
}

\newcommand{\SEPATTERN}{
    \draw[step=1, xshift=14pt, yshift=14pt, \cfill, line cap=round] (0,0) grid (3,3);
    \draw[step=1, xshift=14pt, yshift=14pt, thick] (0,1) -- (3,1);
    \draw[step=1, xshift=14pt, yshift=14pt, thick] (1,0) -- (1,3);
    \foreach \x/\y in {1/2,2/3,3/1} \node[disc, fill=black] at (\x,\y) {};
}

\newcommand{\pattern}{\!\raisebox{-0.3em}{
  \begin{tikzpicture}[line width=0.7pt, scale=0.15]
    \tikzstyle{disc} = [circle,thin,draw=black, minimum size=1.7pt, inner sep=0pt ]
    \PATTERN
  \end{tikzpicture}}
}

\newcommand{\sepattern}{\!\raisebox{-0.3em}{
  \begin{tikzpicture}[line width=0.7pt, scale=0.15]
    \tikzstyle{disc} = [circle,thin,draw=black, minimum size=1.7pt, inner sep=0pt ]
    \SEPATTERN
  \end{tikzpicture}}
}

\newcommand{\subi}{\textsf{Rem1}}
\newcommand{\subii}{\textsf{Rem2}}
\newcommand{\subiii}{\textsf{Rem3}}

\newcommand{\addi}{\textsf{Add1}}
\newcommand{\addii}{\textsf{Add2}}
\newcommand{\addiii}{\textsf{Add3}}

\newcommand{\cfill}{black!40}
\newcommand{\cfilll}{white}

\newcommand{\ns}{4pt}
\newcommand{\nodestyle}{\tikzstyle{every node} = [font=\footnotesize]}
\newcommand{\discstyle}{\tikzstyle{disc} =
  [ circle,thin,fill=\cfilll,draw=black, minimum size=\ns, inner sep=0pt ] }
\newcommand{\style}{
  \nodestyle
  \discstyle
}

\begin{document}

\begin{center}
{\large\bf Equidistributed statistics on   Fishburn matrices and  permutations}
\end{center}

\begin{center}
{\small Dandan Chen$^\dag$,  Sherry H.F. Yan$^\dag$,   Robin D.P. Zhou$^\ddag$}

$^{ \dag}$Department of Mathematics\\
Zhejiang Normal University\\
 Jinhua 321004, P.R. China

$^\ddag$College of Mathematics Physics and Information\\
Shaoxing University\\
Shaoxing 312000, P.R. China

huifangyan@hotmail.com

\end{center}

\begin{abstract}
Recently, Jel\'inek conjectured that there exists a bijection between  certain restricted permutations and Fishburn matrices such that the bijection verifies the    equidistribution   of several statistics.  The main objective of this paper is to establish such a bijection.

\end{abstract}

\noindent {\bf Keywords}: ascent sequence, pattern avoiding permutation,  Fishburn matrix.

\noindent {\bf AMS  Subject Classifications}:  05A15, 05A17, 06A07

\section{Introduction}
Given a sequence of integers   $x=x_1x_2\cdots x_n$, we say
that the sequence $x$ has an {\em ascent} at position $i$ if
$x_i<x_{i+1}$. Let $ASC(x)$ denote the set of the ascent positions of $x$  and let  $asc(x)$ denote the number of ascent of $x$.
A sequence $x=x_1x_2\cdots x_n$ is said to be an {\em ascent
sequence of length $n$} if it satisfies $x_1=0$ and $0\leq x_i\leq
asc(x_1x_2\cdots x_{i-1})+1$ for all $2\leq i\leq n$.
 Let $\mathcal{A}_n$ be the set of ascent sequences of length $n$.
For example,
\[
\mathcal{A}_3=\{000,001,010,011,012\}
\]

Ascent sequences were introduced by Bousquet-M\'elou et al. \cite{Bousquet}
to unify three other combinatorial structures: $(2+2)$-free posets, a family of
permutations avoiding a certain pattern and a class of involutions introduced
by Stoimenow \cite{Stoimenow}.
To be specific, Bousquet-M\'elou et al. \cite{Bousquet} constructed a bijection
between ascents sequences and pattern avoiding permutations, a bijection between
ascent sequences and $(2+2)$-free posets and a bijection between $(2+2)$-free posets
and Stoimenow's involutions.
Dukes and Parviainen \cite{Dukes} completed the results of \cite{Bousquet} by constructing a bijection between ascent sequences and Fishburn matrices.
Hence, all these combinatorial objects are enumerated by  the  Fishburn number $F_n$ (sequence A022493 in OEIS \cite{sloane} )
for memory of Fishburn's
pioneering  work on the interval orders \cite{Fishburn1,Fishburn2, Fishburn3}.
More examples of Fishburn objects  are constantly being discovered.
   Levande \cite{ Levande} introduced the notion of Fishburn diagrams and proved that Fishburn diagrams are counted by Fishburn numbers, confirming a conjecture posed by Claesson and  Linusson \cite{Claesseon}.
Jel\'inek \cite{Jelinek} showed  that some Fishburn triples are
enumerated by Fishburn numbers.

Zagier \cite{Zagier} and Bousquet-M\'elou et al. \cite{Bousquet} obtained the generating function of $F_n$, that is
\[
\sum_{n\geq 0}F_nx^n=\sum_{n \geq 0}\Pi_{k=1}^n(1-(1-x)^k).
\]
Kitaev and Remmel \cite{Kitaev} extended the work and found
the generating function for $(2+2)$-free posets when four statistics are
taken into account.
 Levande \cite{ Levande} and Yan \cite{Yan} independently  presented a combinatorial proof of a conjecture of
Kitaev and Remmel \cite{Kitaev} concerning the generating function for the number of $(2+2)$-free posets.

Let us recall the notions of pattern avoiding permutations and Fishburn matrices
before we state our main results.
Let $S_n$ be the symmetric group on $n$ elements and $\pi=\pi_1\pi_2\cdots\pi_n$ be a permutation of $S_n$.
We say that $\pi$ contains the pattern \mbox{\pattern} if there is a subsequence
$\pi_i\pi_{i+1}\pi_j$  of $\pi$ satisfying that $\pi_i+1=\pi_j<\pi_{i+1}$,
otherwise we say that $\pi$ avoids the pattern  \mbox{\pattern}.
For example, the permutation $42513$ contains the pattern  \mbox{\pattern}
while the permutation $52314$ avoids it.

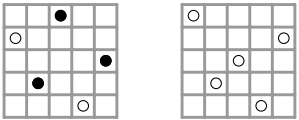
\begin{tikzpicture}[line width=1.1pt, scale=0.3]
  \style
  \draw[step=1, xshift=14pt, yshift=14pt, \cfill, line cap=round]
  (0,0) grid (5,5);
  \foreach \x/\y in {1/4,2/2,3/5,4/1,5/3}
  \node[disc] (\y) at (\x,\y) {};
  \foreach \y in {3,5,2}
  \node[disc,fill=black] at (\y) {};
\end{tikzpicture}\quad\quad
\begin{tikzpicture}[line width=1.1pt, scale=0.3]
  \style
  \draw[step=1, xshift=14pt, yshift=14pt, \cfill, line cap=round]
  (0,0) grid (5,5);
  \foreach \x/\y in {1/5,2/2,3/3,4/1,5/4}
 \node[disc] (\y) at (\x,\y) {};
\end{tikzpicture}

The pattern \mbox{\sepattern} can be defined similarly.
Let $S_n(\pattern)$ be the set of
$(\pattern)$-avoiding permutations of $[n]$ and
 $S_n(\sepattern)$ be the set of
$(\sepattern)$-avoiding permutations of $[n]$, respectively.
These two sets are both enumerated by Fishburn numbers \cite{Bousquet,Parviainen}.
In a permutation $\pi$, we say  $\pi_i$ is a left-to-right maximum (or LR-maximum) if $\pi_i$ is larger than any element among $\pi_1, \pi_2,\ldots,\pi_{i-1}$.
Let $LRMAX(\pi)$ denote the set of LR-maxima of $\pi$ and let
 $LRmax(\pi)$ denote the number of LR-maxima of $\pi$.   Analogously, we can define LR-minima, RL-maxima, RL-minima of a permutation $\pi$.
Denote by $LRMIN(\pi$),  $RLMAX(\pi)$ and  $RLMIN(\pi)$  the set of LR-minima, RL-maxima and  RL-minima of $\pi$, their cardinalities  being denoted by
  $LRmin (\pi)$,  $RLmax(\pi)$ and  $RLmin (\pi)$, respectively.

Fishburn matrices were introduced by Fishburn \cite{Fishburn3} to represent
interval orders.
A Fishburn matrix is an upper triangular matrix with nonnegative integers whose every row and every column contain at least one non-zero entry.
The weight of a matrix is the sum of its entries.
Similarly, the weight of a row (or a column) of a matrix is the sum of the entries in this row (or column).
Denote by $\mathcal{M}_n$   the set of Fishburn matrices of weight $n$.
For example,
\[
\mathcal{M}_3=\{
\begin{pmatrix} 3 \end{pmatrix},
\begin{pmatrix} 2 & 0 \\0& 1 \end{pmatrix},
\begin{pmatrix} 1 & 1 \\0& 1 \end{pmatrix},
\begin{pmatrix} 1 & 0 \\0& 2 \end{pmatrix},
\begin{pmatrix} 1 & 0&0 \\0& 1&0\\0&0&1 \end{pmatrix}
\}.
\]

Given a matrix $A$, we use the term  {\em cell $(i,j)$} of $A$ to refer to the the entry in the $i$-th row and $j$-th column of $A$, and we let    $A_{i,j}$ denote its value. We assume that the rows  of a matrix are numbered from top to bottom  and the columns of a matrix are numbered from left to right in which the topmost row is numbered by $1$ and the leftmost column is numbered by $1$.  A cell $(i,j)$ of a matrix $A$  is said to be zero if $A_{i,j}=0$. Otherwise, it is said to be {\em nonzero}.  A row ( or column) is said be zero if it contains no nonzero cells. Otherwise, it is  said to be {\em nonzero} row ( or  column).

    A cell $(i,j)$ of a matrix $A$ is
a weakly north-east cell (or wNE-cell) if it  is a nonzero cell and any other cell weakly north-east form $c$ is a zero cell. More precisely,  a cell $(i,j)$ of a matrix $A$ is
    a wNE-cell if $A_{s,t}=0$ for all $s\leq i$ and $t\geq j$.

Jel\'inek \cite{Jelinek} posed the following conjecture.

\begin{conj}{ \upshape   (See \cite{Jelinek}, Conjecture  4.1)}\label{conj1}
For every $n$, there is a bijection $\alpha$ between $S_n(\pattern)$ and $\mathcal{M}_n$ satisfying that:
\begin{itemize}
  \item LRmax$(\pi)$ is the weight of the first row of $\alpha(\pi)$,
  \item RLmin$(\pi)$ is the weight of the last column of $\alpha(\pi)$,
  \item RLmax$(\pi)$ is the number of wNE-cells of $\alpha(\pi)$,
  \item LRmin$(\pi)$ is the number of nonzero cells of $\alpha(\pi)$ belonging to
        the main diagonal, and
  \item $\alpha(\pi^{-1})$ is obtained from $\alpha(\pi)$ by transposing along the North-East diagonal.
\end{itemize}
\end{conj}

By using generating functions, Jel\'inek \cite{Jelinek} proved the following symmetric joint distribution on $\mathcal{M}_n$.
\begin{thm}{ \upshape   (See \cite{Jelinek},  Theorem $3.7$)}\label{thm1.1}
For any $n$,  the number of wNE-cells and the weight of the first row   have symmetric joint distribution on $\mathcal{M}_n$.
\end{thm}

Jel\'inek \cite{Jelinek} also posed the following weaker conjecture
which can be followed directly from
 Theorem \ref{thm1.1} and  Conjecture \ref{conj1}.

\begin{conj}{ \upshape   (See \cite{Jelinek}, Conjecture  4.2)}\label{conj2}
For any $n$, LRmax and RLmax have symmetric joint distribution on $S_n(\pattern)$.
\end{conj}

The main objective of this paper is   to establish   a bijection   between  $S_n(\pattern)$ and $\mathcal{M}_n$ which
satisfies the former four items of Conjecture \ref{conj1},
thereby confirming  Conjecture \ref{conj2}.

\section{Bijection between  permutations and ascent sequences  }

In this section, we shall construct a bijection $\theta$
between $S_n(\pattern)$ and $\mathcal{A}_n$,  and show that the map $\theta$ proves the equidistribution of two 4-tuples of statistics.

Let $\pi$ be a permutation in $S_n(\pattern)$ and let
$\tau$ be the permutation obtained by  deleting $n$ from $\pi$.
Then we have that $\tau$ is also a permutation in  $S_n(\pattern)$.
If not, we assume that $\tau_i\tau_{i+1}\tau_j$ is a \pattern pattern
in $\tau$.
Since $\pi$ is (\pattern)-avoiding, we have $\pi_{i+1}=n$.
Then $\pi_i\pi_{i+1}\pi_{j+1}$ forms a \pattern pattern in $\pi$,
a contradiction.
This property allows us to construct the permutation of
$S_n(\pattern)$ inductively, starting from the empty permutation
and adding  a new maximal value at each step.

Let $\tau$ be a permutation in $S_{n-1}(\pattern)$.
The positions where we can insert the element $n$ into $\tau$
to obtain a \pattern-avoiding permutation are called
active sites.
The site after the maximal entry $n$ in $\pi$ is always an active site.
We label the active sites in $\pi$ from right to left with $0,1,2$ and
so on.

The bijection $\theta$ between $S_n(\pattern)$ and $\mathcal{A}_n$
can be defined recursively.
Set $\theta(1)=0$.
Suppose that $\pi$ is a permutation in $S_n(\pattern)$
which is obtained from  $\tau$ by inserting the element $n$ into the
$x_n$-th active site of $\tau$.
Then we set $\theta(\pi)= x_1 x_2 \cdots  x_{n-1} x_n $,
where $\theta(\tau)= x_1 x_2 \ldots x_{n-1}$.

\begin{example}
The permutation $85231647$ corresponds to the sequence $0 1 1 0 2 1 0 3$ since it  is obtained  by the following insertion,
where the subscripts indicate the labels of the active sites.
 \begin{align*}
    _1 1 _0
    &\,\xrightarrow{x_2=1}\,  {_2} 2 {_1} 1 {_0} \\
    &\,\xrightarrow{x_3=1}\,  {_2} 2 \;3{_1} 1 _0 \\
    &\,\xrightarrow{x_4=0}\,  {_2} 2\; 3\; 1 {_1} 4 {_0} \\
    &\,\xrightarrow{x_5=2}\,  {_3} 5 {_2} 2\; 3 \; 1 {_1}4{_0} \\
    &\,\xrightarrow{x_6=1}\,  {_3} 5\; 2 \; 3\; 1 {_2} 6{_1}4 {_0} \\
    &\,\xrightarrow{x_7=0}\,  {_3}5\; 2\; 3 \; 1 {_2}6\; 4{_1}7 {_0} \\
    &\,\xrightarrow{x_8=3}\,  {_4}8{_3}5\; 2\; 3 \; 1 {_2}6\; 4{_1}7 {_0} .
  \end{align*}
\end{example}

\begin{lem}\label{lem:1}
Let $\pi=\pi_1\pi_2\cdots\pi_n$ be a permutation in $S_n(\pattern)$ and
$\theta(\pi)=x= x_1 x_2 \cdots x_n $.
Then we have that
\begin{equation}\label{equ:1}
s(\pi)=2+asc(x)\quad \text{and} \quad a(\pi)=x_n,
\end{equation}
where $s(\pi)$ denotes the number of active sites of $\pi$
and $a(\pi)$ denotes the label of the site located just after the entry $n$ of $\pi$.
\end{lem}

\pf
Suppose that $\pi$ is obtained from $\tau$ by inserting the element $n$ into the
$x_n$-th active site of $\tau$.
Then we have $\theta(\tau)=x'$, where
$x'=x_1x_2\cdots x_{n-1}$.
For any entry $i$ which is to the right of $n$,
$i$ is followed by an active site in $\pi$
if and only if $i$ is followed by an active site in $\tau$.
Since the site after $n$ in $\pi$ is always active, we obtain $a(\pi)=x_n$

Now let us focus on the equation $s(\pi)=2+asc(x)$.
We will prove it by induction on $n$.
It obviously hold for $n=1$.
Assume that it holds for $n-1$.
For any entry $i<n-1$, $i$ is followed by an active site in $\pi$
if and only if $i$ is followed by an active site in $\tau$.
The site after $n$ in $\pi$ is always an active site.
Thus, to determine $s(\pi)$, the only question is whether the site after
$n-1$ is active.
We need consider two cases.

\noindent Case 1: If $0\leq x_n\leq a(\tau)=x_{n-1}$, then the entry $n$
in $\pi$ is to the right of $n-1$.
It follows that the site after $n-1$ is not an active cite in $\pi$.
Since the site after $n-1$ is an active cite in $\tau$,
we have that $s(\pi)=s(\tau)$.
By the induction hypothesis, $s(\tau)=2+asc(x')=2+asc(x)$.
Hence we deduce that $s(\pi)=2+asc(x)$.

\noindent Case 2: If  $x_n> a(\tau)=x_{n-1}$, then the entry $n$
in $\pi$ is to the left of $n-1$.
It yields that the site after $n-1$ is also an active cite in $\pi$.
Hence $s(\pi)=s(\tau)+1$.
Since $x_n>x_{n-1}$, we have that $asc(x)=asc(x')+1$.
By the induction hypothesis, $s(\tau)=2+asc(x')$.
Thus we have $s(\pi)=2+asc(x)$.
This completes the proof.
\qed

\begin{thm}\label{thm2.1}
The map $\theta$ is a bijection between $S_n(\pattern)$
and $\mathcal{A}_n$.
\end{thm}

\pf
We prove this conclusion by induction on $n$.
It obviously holds for $n=1$.
Assume that $\theta$ is a bijection between $S_{n-1}(\pattern)$
and $\mathcal{A}_{n-1}$.

We first show that $\theta$ is a map from $S_n(\pattern)$
to $\mathcal{A}_n$.
Let $\pi=\pi_1\pi_2\cdots\pi_n$
be a permutation in $S_n(\pattern)$
which is obtained from $\tau$ by inserting a maximal entry $n$ in
the active site labeled by $x_n$ in $\tau$.
Then $\theta(\pi)=x= x_1 x_2\cdots x_n $,
where $\theta(\tau)=x'= x_1 x_2 \cdots x_{n-1} $.
To prove that $x\in \mathcal{A}_n$,
it suffices to show that $x_n\leq asc(x')+1$.
Recall that the rightmost active site is labeled $0$. Hence the leftmost
active site in $\tau$ is labeled   $s(\tau)-1$.
By the recursive description of the map $\theta$, we have that
$x_n\leq s(\tau)-1$.
From Lemma \ref{lem:1} we see that $s(\tau)=2+asc(x')$.
Thus we have $x_n\leq asc(x')+1$.
Since $x$ encodes the construction of $\pi$, $\theta$
is an injective map from $S_n(\pattern)$
to $\mathcal{A}_n$.

It remains to show that $\theta$ is surjection.
Let $y= y_1 y_2 \cdots y_n $ be an ascent sequence and $p=p_1p_2\cdots p_{n-1}=\theta^{-1}(y')$, where $y'= y_1 y_2 \cdots y_{n-1} $.
From the definition of ascent sequence and Lemma \ref{lem:1},
we have that $y_n\leq asc(y')+1=s(p)-1$.
Let $q$ be the permutation obtained from $p$ by inserting the maximal
entry $n$ into the active site labeled $y_n$ in $p$.
By the construction of the map $\theta$, it can be easily seen that
$\theta(q)=y$.
This concludes the proof.
\qed

Let $x= x_1 x_2 \cdots x_n $ be an ascent sequence in  $\mathcal{A}_n$.
The {\em  modified ascent sequence} of $x$, denoted by $\hat{x}$,  is defined
by the following procedure:\\
for $i\in ASC(x)$\\
\hspace*{0.5cm}  for $j=1,2,\ldots, i-1$\\
 \hspace*{1cm} if $x_j\geq  x_{i+1}$ then $x_j:=x_j+1$.\\
 For example, for $x=01012213$, we have $ASC(x)=\{1,3,4,7\}$ and $\hat{x}=04012213$.
Modified ascent sequence were introduced by  Bousquet-M\'elou et al.,
see more details in \cite{Bousquet}.

For a permutation $\pi=\pi_1\pi_2\cdots\pi_n\in S_n(\pattern)$, let
$l(\pi_i)$ be the largest label of the active site to the right of $\pi_i$ and
let $LMAXL(\pi)$ be the multiset of $l(\pi_i)$ when $\pi_i$ ranges over all
LR-maxima of $\pi$.
That is $$LMAXL(\pi)=\{l(\pi_i)\mid \pi_i\in LRMAX(\pi)\}.$$
Similarly, let
$$RMAXL(\pi)=\{l(\pi_i)\mid \pi_i\in RLMAX(\pi)\}.$$
Define $$\delta(\pi,q)=\sum_{i\in LMAXL(\pi)}q^i.$$
For example, for $\pi=42178536$, its active sites are labelled as $ {_4}421_378_253_16_0$. Then we have $RMAXL(\pi)=\{0,2\}$ and $LMAXL(\pi)=\{2,2,3\}$.

For an ascent sequence $x= x_1 x_2 \cdots x_n $,
let $zero(x)$ denote the number of zeros in  $x$ and
let $max(x)$ denote the number of elements $x_i$ satisfying
$x_i=asc(x_1x_2\cdots x_{i-1}) +1$.

For a sequence $x= x_1 x_2 \cdots  x_n $,  let
 $$RMIN(x)=\{x_i\mid x_i<x_j \,\, \mbox{for all }  j>i\},$$

$$RMAX(x)=\{x_i\mid x_i\geq x_j \,\, \mbox{for all }  j>i\},$$
and

$$ \chi( x, q)=\sum_{x_i \in RMAX(x)}q^{x_i} .$$

 Denote by $Rmin(x)$ and $Rmax(x)$ the cardinalities of  the set $RMIN(x)$ and $RMAX(x)$, respectively.
\begin{thm} \label{thm2.2}
For any $\pi=\pi_1\pi_2\cdots\pi_n\in S_n(\pattern)$ and $x= x_1 x_2 \cdots x_n \in \mathcal{A}_n$ with $\theta(\pi)=x$, we have
\begin{enumerate}
  \item[(1)] $RLmin(\pi)=zero(x)$;
  \item [(2)] $LRmin(\pi)=max(x)$;
  \item [(3)] $RMAXL(\pi)=RMIN(x)$;
  \item [(4)] $\delta(\pi,q)=\chi(\hat{x},q)$;
  \item [(5)] $RLmax(\pi)=Rmin(x)$;
  \item [(6)] $LRmax(\pi)= Rmax(\hat{x})$.
\end{enumerate}
\end{thm}

\pf Point (5) follows directly from point (3).
Similarly, point (6) is an immediate consequence of the point (4) with $q=1$.
Now we will prove point (1)-(4) by induction on $n$.
It is easily checked that the statement holds for $n=1$.
Assume that it also holds for  some $n-1$ with $n\geq 2$.
Let $\tau$ be the permutation which is obtained from $\pi$ by deleting the
largest entry $n$ in $\pi$.
Then we have that $x'= x_1 x_2 \cdots x_{n-1} =\theta(\tau)$.
From the construction of the bijection $\theta$ and the induction hypothesis,
one can easily verify that

$$
RLmin(\pi)=\left\{
\begin{array}{ll}
RLmin(\tau)+1=zero(x')+1=zero(x)  & \,\, \mbox{if } x_n=0,\\
  RLmin(\tau)=zero(x')=zero(x)&   \,\,  \mbox{otherwise }, \\
\end{array}
\right.
$$

$$
LRmin(\pi)=\left\{
\begin{array}{ll}
LRmin(\tau)=max(x')=max(x)  & \,\,   \mbox{if } x_n\leq asc(x'),   \\
  LRmin(\tau)+1=max(x')+1=max(x)&   \,\,  \mbox{if } x_n=asc(x')+1, \\
\end{array}
\right.
$$

and
$$
\begin{array}{lll}
RMAXL(\pi) &=&
 \{i\mid i\in RMAXL(\tau), i<x_n\}\cup\{x_n\}\\
 &=&\{i\mid i\in RMIN(x'), i<x_n\}\cup\{x_n\}\\
 &=&RMIN(x).
 \end{array}
 $$

 For  point (4), we consider two cases.
 If $x_n\leq x_{n-1}$, then $n$ is to the right of $n-1$ in $\pi$.
  Notice that all the LR-maxima in $\tau$ are also LR-maxima in $\pi$. One can easily check that
  $LMAXL(\pi)=LMAXL(\tau)\cup \{x_n\}$ and    $RMAX(\hat{ x })=RMAX(\hat{x'})\cup \{x_n\}$.
 Hence we have
 $$
 \delta(\pi,q)=\delta(\tau,q)+q^{x_n}=\chi(\hat{x'},q)+q^{x_n}=\chi(\hat{x},q).
 $$
 If $ x_n> x_{n-1}$,  then $n$ is to the left of $n-1$ in $\pi$.
 In this case, $\tau_i$ is a LR-maximum in $\pi$ if and only if
 $\tau_i$ is a LR-maximum in $\tau$ and $l(\tau_i)\geq x_n$.
 After the inserting $n$ into $\tau$, $l(\tau_i)$ is increased by $1$ if
 $\tau_i$ is also a LR-maximum in $\pi$.
 Hence we have that
 $$
 \delta(\pi,q)=\sum_{i\in LMAXL(\tau),i\geq x_n}q^{i+1}+q^{x_n}=
 \sum_{i\in RMAX(\hat{x'}),i\geq x_n}q^{i+1}+q^{x_n}=
 \chi(\hat{x},q),
 $$
 where the last equality follows from the fact that $$RMAX(\hat{x})=\{i+1\mid i\in RMAX(\hat{x'}), i\geq x_n\}\cup \{x_n\}.$$
 This completes the proof. \qed

Combining  Theorems \ref{thm2.1} and \ref{thm2.2}, we are led to the following result.

\begin{thm}\label{thm2.3}
The map $\theta$ is a bijection between $S_n(\pattern)$ and $\mathcal{A}_n$.  Moreover, for any $\pi \in S_n(\pattern)$ and $x\in \mathcal{A}_n$ with $\theta(\pi)=x$, we have
$$
(RLmin, LRmin, RLmax)\pi= (zero, max, Rmin)x
$$
and $LRmax(\pi)=Rmax(\hat{x}) $.
  \end{thm}

\section{Bijection  between
 ascent sequences and Fishburn matrices}

The main objective of this section is to establish a bijection $\phi$ between
$\mathcal{A}_n$ and $\mathcal{M}_n$.
To this end,    we will define a removal operation and an
addition operation on the matrices of $\mathcal{M}_n$.

Given a matrix $A$ in $\mathcal{M}_n$, let $dim(A)$ denote the number of rows of
the matrix $A$ and let $index(A)$ denote the smallest value of $i$ such that $A_{i,dim(A)}>0$. Denote by $rsum_i(A)$ and $csum_i(A)$  the sum of the entries in row $i$ and column
$i$ of $A$, respectively.
We define a removal operation $f$ on a given matrix $A\in \mathcal{M}_n$ as follows.

\begin{enumerate}
\item[(\subi)] If $rsum_{index(A)}(A)>1$,  then let $f(A)$ be the matrix $A$ with the entry $A_{index(A),dim(A)}$ reduced by $1$.
\item[(\subii)] If $rsum_{index(A)}(A)=1$ and $index(A)=dim(A)$, then let $f(A)$ be the matrix $A$ with row $dim(A)$ and column $dim(A)$ removed.
\item[(\subiii)]
If $rsum_{index(A)}(A)=1$ and $index(A)<dim(A)$, then we construct $f(A)$ in the following way. Let $S$ be the set of indices $j$ such that $j\geq index(A)$    and column $j$ contains  at least one nonzero entry   above row $index(A)$. Suppose that $S=\{c_1, c_2, \ldots, c_\ell\}$ with $c_1<c_2\ldots<c_\ell$. Clearly we have $c_1=index(A)$.  Let $c_{\ell+1}=dim(A)$.   For all $1\leq i< index(A)$ and $1\leq j\leq \ell$, move all the entries in the cell $(i, c_{j})$ to the cell $(i, c_{j+1})$.
  Simultaneously  delete row $index(A)$ and column $index(A)$.

\end{enumerate}

\begin{example}
Let $A,B,C$ be the following three Fishburn matrices:
$$
A=\left(
  \begin{array}{cccc}
    1 & 2 & 0 & 0 \\
    0 & 2 & 1 & 0 \\
    0 & 0 & 2 & 1 \\
    0 & 0 & 0 & 2 \\
  \end{array}
\right);
\quad
B=\left(
  \begin{array}{cccc}
    1 & 0 & 2 & 0 \\
    0 & 3 & 0 & 0 \\
    0 & 0 & 2 & 0 \\
    0 & 0 & 0 & 1 \\
  \end{array}
\right);
\quad
C= \left(
            \begin{array}{ccccc}
              2 & 4 & 1 & 3 &  0  \\
              0 & 5 & 2 & 2 &  0  \\
              0  &  0  &  0  &  0  & 1  \\
              0 & 0 &  0  & 1 & 3 \\
              0 & 0 &  0  & 0 & 2 \\
            \end{array}
          \right)
.$$
\end{example}

For Matrix $A$, rule (\subi) is applied since $rsum_{index(A)}(A)=3$ and
$$f(A)=\left(
  \begin{array}{cccc}
    1 & 2 & 0 & 0 \\
    0 & 2 & 1 & 0 \\
    0 & 0 & 2 & 0 \\
    0 & 0 & 0 & 2 \\
  \end{array}
\right).$$

For Matrix $B$,
since $rsum_{index(B)}(B)=1$ and $index(B)=dim(B)$, rule (\subii) is applied
and
$$f(B)=\left(
  \begin{array}{ccc}
    1 & 0 & 2  \\
    0 & 3 & 0  \\
    0 & 0 & 2  \\
      \end{array}
\right).$$

For matrix $C$, since $rsum_{index(C)}(C)=1$ and $index(C)<dim(C)$,
rule (\subiii) is applied. It is easy to check that $S=\{3,4\}$,  and thus we have
$$
f(C)=\left(
            \begin{array}{cccc}
              2 & 4 & 1 & 3  \\
              0 & 5 & 2 & 2  \\
              0 & 0 & 1 & 3 \\
              0 & 0 & 0 & 2 \\
            \end{array}
          \right).
$$

The following lemma shows that the removal operation on a Fishburn matrix of $\mathcal{M}_n$ will yield a Fishburn matrix in $\mathcal{M}_{n-1}$.

\begin{lem}\label{lem1}
Let $n\geq 2$ be an integer and $A\in \mathcal{M}_n$, then we have
that $f(A)\in \mathcal{M}_{n-1}$.
\end{lem}

\pf
It is easily seen that for any removal operation applied on the matrix $A$,
 the weight of $f(A)$ is one less than the weight of $A$.
 It is trivial to check that there exists no  zero columns or rows  in $f(A)$. Moreover, the removal operation also preserves the property of being   upper-triangular.
 Thus, $f(A) \in \mathcal{M}_{n-1}$.
 This completes the proof.
 \qed

 Lemma \ref{lem1} tells us  that for any $A\in \mathcal{M}_n$, after $n$ applications of the removal operation  $f$ to $A$, we will get   a sequence of Fishburn  matrices, say $A^{(1)}, A^{(2)}, \ldots,A^{(n)}$, where $A^{(k-1)}=f(A^{(k)})$ for all $1<k\leq n$ and $A^{(n)}=A$. Define $\psi(A)=x=x_1x_2\ldots x_n$  where
 $x_k=index(A^{(k)})$.

We now define an addition operation $g$ on a Fishburn matrix which
   is  shown to be the  inverse of the removal operation later.
  Given a matrix $A\in \mathcal{M}_n$ and $i\in [0, dim(A)]$,
 We construct a matrix $g(A, i)$  in the following manner.

 \begin{enumerate}
   \item[(\addi)] If $0\leq i\leq index(A)-1$, then let $g(A, i)$ be the matrix obtained from $A$ by increasing the entry in the cell $(i+1, dim(A))$ by 1.
   \item[(\addii)]If $i=dim(A)$, then let $g(A, i)$ be the matrix $\left(
            \begin{array}{cc}
              A & 0  \\
              0 & 1  \\
            \end{array}
          \right). $
   \item[(\addiii)] If $index(A)\leq i<dim(A)$,  then  we construct $g(A,i)$ in the following way. In $A$, insert a new (empty) row between rows $i$ and $i+1$, and insert a new (empty) column between columns $i$ and $i+1$. Let the new row be filled with all zeros  except for the rightmost cell which is filled with a $1$.  Denote by $A'$ the resulting matrix.  Let $T$ be the set of indices $j$ such that $j\geq i+1$ and  column $j$  contains  at least  one nonzero cell above row $i+1$. Suppose that $T=\{c_1, c_2, \ldots, c_\ell\}$. Clearly we have $c_{\ell}=dim(A')$. Let $c_0=i+1$.  For all $1\leq a\leq i$ and $1\leq b\leq \ell $, move all the entries in the cell $(a, c_{b})$ to the cell $(a, c_{b-1})$, and fill  all the cells  which are in column $ dim(A') $  and  above row $i+1$ with zeros.

 \end{enumerate}

\begin{example}
Consider the matrix
$$
A=\left(
            \begin{array}{cccc}
              2 & 4 & 0 & 3  \\
              0 & 5 & 0 & 2  \\
              0 & 0 & 1 & 3 \\
              0 & 0 & 0 & 2 \\
            \end{array}
          \right).
$$
Obviously, we have $dim(A)=4$ and $index(A)=1$.
For $i=0$, since $i\leq index(A)-1$, rule (\addi) applies and we get
$$
g(A,0)=\left(
            \begin{array}{cccc}
              2 & 4 & 0 & 4  \\
              0 & 5 & 0 & 2  \\
              0 & 0 & 1 & 3 \\
              0 & 0 & 0 & 2 \\
            \end{array}
          \right).
$$
For $i=4$, since $ i=dim(A)$, rule (\addii)  applies and we get
$$
g(A,4)=\left(
            \begin{array}{ccccc}
              2 & 4 & 0 &3 & 0 \\
              0 & 5 & 0 & 2& 0  \\
              0 & 0 & 1 & 3 & 0\\
              0 & 0 & 0 & 2 & 0\\
               0 & 0 & 0 & 0 & 1\\
            \end{array}
          \right).
$$
For $i=1$,    since $index(A) \leq i<dim(A)$,    rule (\addiii)  applies and we get
$$
A'=\left(
            \begin{array}{ccccc}
              2 & \bf{0} & 4 & 0 & 3 \\
             \bf{0} &\bf{0} & \bf{0} & \bf{0}& \bf{1}  \\
              0 & \bf{0} & 5 & 0 &2\\
              0 & \bf{0} & 0 & 1 & 3\\
               0 & \bf{0} & 0 & 0 & 2\\
            \end{array}
          \right),
$$
where the new inserted row and column are illustrated in bold.
Then we have $T=\{3,5\}$. Finally, we get
$$
g(A,1)=\left(
            \begin{array}{ccccc}
              2 & 4 & 3 & 0 & 0 \\
              0 & 0 & 0 & 0& 1  \\
              0 & 0 & 5 & 0 &2\\
              0 & 0 & 0 & 1 & 3\\
               0 & 0 & 0 & 0 & 2\\
            \end{array}
          \right).
$$
\end{example}

By similar arguments as in the  proof of Lemma \ref{lem1}, one can easily verify that the addition  operation  will also yield  a Fishburn matrix.

\begin{lem}\label{lem2}
For any   matrix $A\in \mathcal{M}_{n-1}$ and $i\in [0, dim(A)]$, we have
that $g(A, i)\in \mathcal{M}_{n}$.
\end{lem}

We now define a map $\phi$ from $\mathcal{A}_n$ to $\mathcal{M}_n$ recursively as follows. Given an ascent sequence $x=x_1  x_2  \ldots, x_n $, we define $A^{ (1) }=(1)$ and $A^{ (k ) }=g(A^{ (k-1) }, x_k)$ for all $1<k\leq n$.  Set $\phi(x)=A^{(n)}$.

Next we aim to show that the map $\phi$  is well defined and has the following desired properties.

\begin{lem}\label{lem3}
For any $x= x_1  x_2  \cdots  x_n \in \mathcal{A}_n$, we have $\phi(x)\in \mathcal{M}_n$ satisfying that  $dim(\phi(x))=asc(x)+1\,\,\, \mbox{and} \,\, index(\phi(x))=x_n+1$.
\end{lem}

\pf  We will prove by induction on $n$.
It is trivial to check that the statement holds for $n=1$. Assume that it also holds for   $n-1$, that is,
   $$
\phi(x')\in \mathcal{M}_{n-1},  \,\, dim(\phi(x'))=asc(x')+1\,\,\, \mbox{and} \,\, index(\phi(x'))=x_{n-1}+1,
$$
where $x'= x_1  x_2  \cdots  x_{n-1} $.
Since $0\leq x_n\leq asc(x')+1=dim(\phi(x'))$,  from Lemma \ref{lem2} we see that $\phi(x)=g(\phi(x'), x_n)\in \mathcal{M}_n $.  From the construction of the addition operation, one can easily verify that $index(\phi(x))=x_{n}+1$ and
$$
dim(\phi(x))=\left\{
\begin{array}{ll}
dim(\phi(x'))=asc(x')+1=asc(x)+1  & \,\, \mbox{if } x_n\leq x_{n-1},\\
dim(\phi(x'))+1=asc(x')+2=asc(x)+1 &   \,\, \mbox{if } x_n> x_{n-1}.\\
\end{array}
\right.
$$
The result follows. \qed

For a matrix $A$, let $ NE(A)=\{i-1| \,\mbox{the cell }\, (i,j) \,\, \mbox{is a wNE-cell of}\, A\}$ and let   $ne(A)$  denote the number of  wNE-cells of $A$.
 Define
$$\lambda(A, q)=\sum_{i=1}^{dim(A)}A_{i,dim(A)}q^{i-1} .$$
Denote by $tr(A)$ the number of nonzero cells belonging to the main diagonal of $A$.

\begin{lem}\label{lem5}
 For any $x= x_1 x_2  \cdots  x_n \in \mathcal{A}_n$ and $A\in \mathcal{M}_n$ with $A=\phi(x)$, we have the following relations.
 \begin{itemize}
 \item[(1)] $zero(x)=rsum_1(A)$;
 \item[(2)] $max(x)=tr(A)$;
 \item[(3)] $ RMIN(x)= NE(A)$;
 \item[(4)]$\chi( \hat{x}, q)=\lambda(A, q)$;
 \item[(5)] $Rmin(x)=ne(A)$;
 \item[(6)] $Rmax(\hat{x})=csum_{dim(A)}(A)$.
   \end{itemize}
\end{lem}

\pf   Point (5) follows directly  from point (3). Similarly,  point (6) is an immediate consequence of the proof  of point (4) with $q=1$. Now we verify points (1)-(4) by induction on $n$. Clearly, the statement holds for $n=1$. Assume that it also holds for any some $n-1$ with $n\geq 2$. Let $x'=x_1x_2\cdots x_{n-1}$ and $B=\phi(x')$.
Recall that $A=g(B, x_n)$.  From the definition of the addition operation $g$ and the induction hypothesis, it is not difficult to verify that
$$
rsum_{1}(A)=\left\{
\begin{array}{ll}
rsum_1(B)+1=zero(x')+1=zero(x),  & \,\, \mbox{if } x_n=0,\\
  rsum_1(B)=zero(x')=zero(x),&   \,\,  \mbox{otherwise }, \\
\end{array}
\right.
$$
and
$$
tr(A)=\left\{
\begin{array}{ll}
tr(B)=max(x')=max(x)  & \,\,   \mbox{if } x_n\leq asc(x'),   \\
  tr(B)+1=max(x')+1=max(x)&   \,\,  \mbox{if } x_n=asc(x')+1. \\
\end{array}
\right.
$$
For point (3), from the construction of the addition operation $g$, we see that the cell  $(x_n+1, dim(A))$ is always a wNE cell. Moreover, there is a  wNE-cell in row $i$ of $A$ if and only if there is a  wNE-cell in row $i$ of $B$ and $i<x_{n}+1$. This yields that
$$
\begin{array}{lll}
NE(A) &=&
 \{i\mid i\in NE(B), i<x_n\}\cup\{x_n\}\\
 &=&\{i\mid i\in RMIN(x'), i<x_n\}\cup\{x_n\}\\
 & =&RMIN(x).
 \end{array}
 $$
 For point (4), we have two cases.

  If $x_n\leq x_{n-1}=index(B)-1$,  then    rule  (\addi ) applies.
  It is trivial to check that
 $$\lambda(A, q)= q^{x_n}+\lambda(B, q)=q^{x_n}+ \chi(\hat{x'}, q)=\chi(\hat{x}, q),
 $$
 where the last equality follows from the fact that $RMAX(\hat{x})=RMAX(\hat{x})\cup \{x_n\}$.

 If $x_n> x_{n-1}=index(B)-1$,  then  either rule (\addii  )  or rule  (\addiii ) applies.
 It is not difficult to verify that
   $$\lambda(A, q)=  q^{x_n} + \sum_{i\geq x_{n}+1}B_{i,dim(B)}q^{i}=q^{x_n} + \sum_{i \in RMAX(\hat{x'}), i\geq x_n}q^{i+1}=\chi(\hat{x}, q),
 $$
 where the last equality follows from the fact that $$RMAX(\hat{x})=\{i+1\mid i \in RMAX(\hat{x'}), i\geq x_n\}\cup  \{x_n\}.$$
     This completes the proof. \qed

 \begin{lem}\label{lem4}
 For any $x= x_1  x_2  \ldots x_n \in \mathcal{A}_n$, we have $\psi(\phi(x))=x$.
 \end{lem}
\pf
Suppose that we get a sequence of matrices $ A^{(1)}, A^{(2)},  \ldots, A^{(n)}$ when we apply the map $\phi$ to $x$, where $A^{(1)}=(1)$ and  $ A^{(k)}= g(A^{(k-1)}, x_k)$ for all $1< k\leq n$.
Similarly, suppose that when we apply the map $\psi$ to $\phi(x)$, we get a sequence $y= y_1  y_2 \ldots  y_n $ and a sequence of matrices $B^{(1)}, B^{(2)}, \ldots, B^{(n)}$, where  $B^{(n)}=\phi(x)$,  $B^{(k)}=f(B^{(k+1)})$ for all $1\leq k<n$,  and $y_{k}=index(B^{(k)})-1$.   Lemma \ref{lem3} ensures that $index(A^{(k)})=x_k+1$. In order to prove $x=y$, it suffices to show that $A^{(k)}=B^{(k)}$  for all $1\leq k\leq n$. We proceed to prove this assertion by induction on $n$. Clearly, we have $B^{(n)}=\phi(x)=A^{(n)}$. Assume that we have $A^{(j)}=B^{(j)}$ for all $j\geq k+1$.  In the following we aim to show that $A^{(k)}=B^{(k)}$. By the induction hypothesis, it suffices to show that $f(A^{(k+1)})=A^{(k)}$.  We have three cases.

Let us assume that $0\leq x_{i+1}<index(A^{(k)})$. Then  rule (\addi) applies and $A^{(k+1)}$ is simply a copy of $A^{(k)}$ with the entry in the cell $(x_{i+1}+1, dim(A^{(k)}))$   increased by one.  Clearly, we have $dim(A^{(k)})=dim(A^{k+1})$,  $index(A^{(k+1)})=x_{i+1}+1$ and  $rsum_{x_{i+1}+1}(A^{(k+1)})>1$.  So rule  (\subi) applies and $f(A^{(k+1)})$ is obtained from $A^{(k+1)}$ by decreasing the the entry in the cell $(x_{i+1}+1, dim(A^{(k+1)}))$ by one. Thus we have $f(A^{(k+1)})=A^{(k)}$.

Next assume that $x_{i+1}=dim(A^{(k)})$. Then rule (\addii) applies and $A^{(k+1)}=\left(
            \begin{array}{cc}
              A^{(k)} & 0  \\
              0 & 1  \\
            \end{array}
          \right). $
In this case,   we have     $index(A^{(k+1)})=x_{i+1}+1=dim(A^{(k+1)})$ and  $rsum_{x_{i+1}+1}(A^{(k+1)})=1$. So rule  (\subii) applies and $f(A^{(k+1)})$ is obtained from $A^{(k+1)}$ by  removing column $dim(A^{(k+1)})$ and row $dim(A^{(k+1)})$. Thus we have $f(A^{(k+1)})=A^{(k)}$.

 If $index(A^{(k)})\leq x_{i+1}<dim(A^{(k)})$, then rule (\addiii) applies and $A^{(k+1)}$ is obtained from $A^{(k)}$ in the following way. First we insert a new (empty) row between rows $x_{i+1}$ and $x_{i+1}+1$, and insert a new (empty) column between columns $x_{i+1}$ and $x_{i+1}+1$. Let the new row be filled with all zeros  except for the rightmost cell which is filled with a $1$.  Denote by $A'$ the resulting matrix.  Let $T$ be the set of indices $j$ such that $j\geq x_{i+1}+1$ and  column $j$  contains  at least  one nonzero cell above row $x_{i+1}+1 $. Suppose that $T=\{c_1, c_2, \ldots, c_\ell\}$ with $c_1<c_2<\ldots<c_{\ell}$.   Let $c_0=x_{i+1}+1 $.  For all $1\leq a\leq x_{i+1}$ and $1\leq b\leq \ell $, move all the entries in the cell $(a, c_{b})$ to the cell $(a, c_{b-1})$, and fill  all the cells in column $ dim(A') $ and above row $x_{i+1}+1$ with zeros. It is easy to check that $dim(A^{(k+1)})=dim(A^{(k)})+1$,  $index(A^{(k+1)})=x_{i+1}+1$ and $rsum_{x_{i+1}+1}(A^{(k+1)})=1$. So rule (\subiii) applies and $f(A^{(k+1)})$ is obtained from $A^{(k+1)}$ by the following procedure. Let $S$ be the set of indices $j$ such that $j\geq x_{i+1}+1$    and column $j$ contains  at least one nonzero entry   above row $x_{i+1}+1$. It is not difficult to check that $S=\{c_0, c_1, c_2, \ldots, c_{\ell-1}\}$.    Let $c_{\ell}=dim(A^{(k+1)})$.   For all $1\leq a< x_{i+1}-1$ and $1\leq b\leq \ell-1$, move all the entries in the cell $(a, c_{b})$ to the cell $(a, c_{b+1})$.
  Simultaneously  delete row $x_{i+1}+1$ and column $x_{i+1}+1$. These operations simply reverse the construction of $A^{(k+1)}$ from $A^{(k)}$,
and therefore $f(A^{(k+1)})=A^{(k)}$.  This completes the proof. \qed

\begin{thm}\label{thm3.1}
The map $\phi$ is a bijection between $\mathcal{A}_n$ and $\mathcal{M}_n$. Moreover, for any $x\in \mathcal{A}_n$ and $A\in \mathcal{M}_n$ with $\phi(x)=A$, we have
$$
(zero, max, Rmin)x=(rsum_1, tr, ne )A
$$
and $Rmax(\hat{x})=csum_{dim(A)}(A)$.
\end{thm}
\pf By Lemma \ref{lem5}, it remains to show that the map $\phi$ is a bijection.   Lemma \ref{lem4} tells us that if  $\phi(x)=\phi(y)$ then we have   $x=y$ for any $x,y\in \mathcal{A}_n$, and thus $\phi$ is injective. And, by cardinality reasons, it follows that $\phi$ is bijective. This completes the proof. \qed

\begin{remark}

  Dukes and Parviainen  \cite{Dukes} defined a bijection $\Gamma$ between $\mathcal{A}_n$ and $\mathcal{M}_n$, and  showed that the bijection $\Gamma$ proves the equidistribution of two triples of statistics, that is, $$
(zero, max)x=(rsum_1, tr  )\Gamma(x)
$$
and $Rmax(\hat{x})=csum_{dim(\Gamma(x))}\Gamma(x)$.  But unlike our bijection $\phi$, the bijection $\Gamma$ does not transform $Rmin$ to $ne$.
\end{remark}

Combining Theorems \ref{thm1.1} and \ref{thm3.1}, we are led to the following symmetric joint distribution on ascent sequences.
\begin{cor}
For any $n$,    the statistics  $zero$ and $Rmin$  have symmetric joint distribution   on $\mathcal{A}_n$.
\end{cor}

Given a  matrix $A\in \mathcal{M}_n $, the {\em flip} of $A$, denoted by $\mathcal{F}(A)$, is the matrix obtained from $A$ by transposing  along the North-East diagonal. It is not difficult to check that for any $A\in \mathcal{M}_n$, we have $\mathcal{F}(A)\in \mathcal{M}_n$ satisfying that
$$
 (rsum_1, tr, ne, csum_{dim(A)} )A= (csum_{dim(\mathcal{F}(A))}, tr, ne, rsum_{1} )\mathcal{F}(A).
$$

In view of Theorems \ref{thm2.3} and \ref{thm3.1}, we are led to the following result,   confirming  the former four items of Conjecture \ref{conj1}.
\begin{thm}
The map $\alpha= \mathcal{ F}\cdot \phi  \cdot \theta $ is a bijection between   $S_n(\pattern)$ and $\mathcal{M}_n$  satisfying that:
\begin{itemize}
  \item LRmax$(\pi)$ is the weight of the first row of $\alpha(\pi)$,
  \item RLmin$(\pi)$ is the weight of the last column of $\alpha(\pi)$,
  \item RLmax$(\pi)$ is the number of wNE-cells of $\alpha(\pi)$,
  \item LRmin$(\pi)$ is the number of nonzero cells of $\alpha(\pi)$ belonging to
        the main diagonal.
\end{itemize}

\end{thm}

\begin{remark}
It should be noted that our bijection $\alpha$ does not  verify the last item of Conjecture \ref{conj1}. For example, let $\pi=85231647$. Then  we have $\pi^{-1}=53472681$, $\theta(\pi)=x=01102103$ and $\theta(\pi^{-1})=y=01223131$. It is easy to check that $asc(x)=3$ and $asc(y)=4$. By Lemma \ref{lem3}, we have $dim(\phi(x))=4$ and $dim(\phi(y))=5$. This implies that the resulting matrices $\alpha(\pi)$ and $\alpha(\pi^{-1})$ have different dimensions, and thus $\alpha(\pi^{-1})\neq \mathcal{ F}(\alpha(\pi))$.
\end{remark}
\vspace{0.5cm}
 \noindent{\bf Acknowledgments.}
This work was supported by the National Science Foundation of China (11671366 and 11626158) and the Zhejiang Provincial Natural Science Foundation of China (   LQ17A010004).

\end{document}